\documentclass[3p,10pt]{elsarticle}

\usepackage{amssymb,latexsym,amsmath,color,subfigure,amsthm}

\journal{}

\newcommand{\eps}{\varepsilon}
\newcommand{\set}[1]{\left\{#1\right\}}
\newcommand{\abs}[1]{\left|#1\right|}
\newcommand{\p}{\partial}
\newcommand{\mB}{\mathbf{B}}
\newcommand{\md}{\mathbf{d}}
\newcommand{\mr}{\mathbf{r}}
\newcommand{\ms}{\mathbf{x}_{\mathrm{s}}}
\newcommand{\mx}{\mathbf{x}}
\newcommand{\vt}{\boldsymbol{\theta}}
\newcommand{\vx}{\boldsymbol{\vartheta}}

\theoremstyle{plain}
\newtheorem{thm}{Theorem}[section]
\newtheorem{lem}[thm]{Lemma}

\theoremstyle{remark}
\newtheorem{rem}{Remark}[section]
\newtheorem{ex}{Example}[section]

\begin{document}

\begin{frontmatter}



\title{Detection of small inhomogeneities via direct sampling method in transverse electric polarization}


\author{Won-Kwang Park}
\ead{parkwk@kookmin.ac.kr}
\address{Department of Information Security, Cryptology, and Mathematics, Kookmin University, Seoul, 02707, Korea}

\begin{abstract}
Various studies have confirmed the possibility of identifying the location of a set of small inhomogeneities via a direct sampling method; however, when their permeability differs from that of the background, their location cannot be satisfactorily identified. However, no theoretical explanation for this phenomenon has been verified. In this study, we demonstrate that the indicator function of the direct sampling method can be expressed by the Bessel function of order one of the first kind and explain why the exact locations of inhomogeneities cannot be identified. Numerical results with noisy data are exhibited to support our examination.
\end{abstract}


\end{frontmatter}





\section{Introduction}
In this study, we consider an inverse-scattering problem for identifying the locations of unknown inhomogeneities with small diameter from measured far-field pattern data. For this purpose, various identification methods have been developed, of which most are based on Newton-type iteration schemes. However, to successfully apply these schemes, one must begin the iteration procedure with a good initial guess that is close to the unknown inhomogeneities. Moreover, it is very difficult to identify multiple inhomogeneities simultaneously using iteration schemes

To quickly identify multiple inhomogeneities, various techniques have been developed; these include MUltiple SIgnal Classification (MUSIC) \cite{AIL1,CZ,P-MUSIC1}, topological derivative \cite{AGJK,BGN,P-TD3}, linear sampling method \cite{CC,HM1,KR}, and Kirchhoff and subspace migrations \cite{AGKPS,KP,P-SUB3}. However, these techniques still require a significant amount of incident-field and corresponding scattered-field directional data to guarantee an acceptable result. By contrast, the direct sampling method needs one or a small number of incident-field data points and has been confirmed to be a very stable and effective detection technique \cite{IJZ1,IJZ2,LZ}. Most research has focused on the detection of targets whose permittivities differ from the background; concerning those whose permeabilities differ from the background, however, little has been performed. Thus, this study analyzes the indicator function of the direct sampling method in a full-view inverse scattering problem. To this end, we construct a relationship using Bessel functions of order one of the first kind. This is based on the fact that a collected far-field pattern can be represented as an asymptotic expansion formula owing to the existence of small inhomogeneities when their permeabilities differ from the background space. From the identified structure, we can examine certain properties of direct sampling method and explain unexplained phenomena.

This paper is organized as follows. In Section \ref{sec:2}, we briefly introduce the two-dimensional direct scattering problem and the asymptotic expansion formula for the far-field pattern. In Section \ref{sec:3}, we establish a relationship between indicator function of direct sampling method and Bessel function of order one and explain its properties. In Section \ref{sec:4}, several results of numerical simulations are presented to support our establishment. A brief conclusion is given in Section \ref{sec:5}.

\section{Direct scattering problem}\label{sec:2}
We briefly discuss two-dimensional time-harmonic electromagnetic scattering from small inhomogeneities located in the homogeneous space $\mathbb{R}^2$. Throughout this study, we assume that all small electromagnetic inclusions $\Sigma_m$, $m=1,2,\cdots,M$, are embedded in $\mathbb{R}^2$ and are characterized by their location $\mx_m$ and size $r_m$:
\[\Sigma_m=\mx_m+r_m\mathbf{B}_m,\]
where $\mathbf{B}_m$ is a simply connected smooth domain containing the origin and $r_m$ is a small positive constant. For the sake, we assume that $\mB_m$ is a unit circle and $r_m$ denotes the radius of $\Sigma_m$. Let $\Sigma$ denote the collection of $\Sigma_m$, and $\omega$ be a given positive angular frequency.

Let us assume that every $\Sigma_m$ is characterized by its magnetic permeability at a given $\omega$ and that all permittivities are equal to $\eps\equiv1$, $\mu_0$ and $\mu_m$ denote the magnetic permeabilities of $\mathbb{R}^2$ and $\Sigma_m$, respectively. With this, we can define the piecewise constant magnetic permeability $\mu(\mx)$ as
\[\mu(\mx)=\left\{\begin{array}{cl}
\mu_0,&\mx\in\mathbb{R}^{2}\backslash\overline{\Sigma}\\
\mu_m,&\mx\in\Sigma_m.
\end{array}\right.\]
We let $k=\omega\sqrt{\eps_0\mu_0}=2\pi/\lambda$ be the wavenumber, where $\lambda$ is the wavelength satisfying $r_m\ll\lambda$ for all $m=1,2,\cdots,M$ (if not, this becomes an imaging of extended target problem, refer to \cite{HSZ1}). Throughout this paper, we assume that $\omega$ is sufficiently large and all $\Sigma_m$ are well-separated such that
\begin{equation}\label{Separated}
k|\mx_m-\mx_{m'}|\gg0.75,
\end{equation}
for $m,m'=1,2,\cdots,M$ and $m\ne m'$.

In this study, we consider the following plane-wave illumination: let $\psi_{\mathrm{inc}}:=\exp(ik\md\cdot\mx)$ be the incident field with direction of propagation $\md\in\mathbb{S}^1$ and $\psi(\mx)$ be the time-harmonic total field that satisfies the Helmholtz equation
\[\nabla\cdot\left(\frac{1}{\mu(\mx)}\nabla\psi(\mx)\right)+\omega^2\psi(\mx)=0\]
with transmission condition on the boundary of $\Sigma_m$. Here $\mathbb{S}^1$ denotes the two-dimensional unit circle. Let $\psi_{\mathrm{scat}}(\mx)$ be the scattered field, which satisfies the Sommerfeld radiation condition
\begin{equation*}
  \lim_{\abs{\mr}\to\infty}\sqrt{\abs{\mr}}\left(\frac{\p\psi_{\mathrm{scat}}(\mx)}{\p\mr}-ik\psi_{\mathrm{scat}}(\mx)\right)=0
\end{equation*}
uniformly in all directions $\mr=\mx/\abs{\mx}$. The far-field pattern of the scattered field satisfies
\begin{equation}\label{FarField}
  \psi_{\mathrm{scat}}(\mx)=\frac{\exp(ik|\mx|)}{\sqrt{|\mx|}}\psi_{\infty}(\md,\vt)+o\left(\frac{1}{\sqrt{|\mx|}}\right)
\end{equation}
uniformly on $\vt=\mx/|\mx|$, as $|\mx|\longrightarrow\infty$.
 As shown in \cite{AK2}, $\psi_{\infty}$ can be written using the following asymptotic expansion formula.

\begin{lem}[Asymptotic formula]
  For sufficiently large $k$, $\psi_{\infty}(\md,\vt)$ can be represented as follows.
  \begin{equation}\label{asymptotic}
    \psi_{\infty}(\md,\vt)\approx-\frac{k^2(1+i)}{4\sqrt{k\pi}}\sum_{m=1}^{M}(r_m)^2|\mB_m|(\md\cdot\mathbb{M}(\mx_m)\cdot\vt)\exp(ik\md\cdot\mx_m)\exp(-ik\vt\cdot\mx_m).
  \end{equation}
  Here, $\mathbb{M}(\mx_m)$ is a $2\times2$ diagonal matrix with components $2\mu_0/(\mu_m+\mu_0)$ and $|\mB_m|=\pi$ is the area of $\mB_m$.
\end{lem}

\section{Indicator function of the direct sampling method}\label{sec:3}
In this section, we introduce an indicator function of the direct sampling method from the collected far-field patterns. Assume that we have the following set of far-field patterns:
\[S:=\set{\psi_{\infty}(\md,\vt_n):n=1,2,\cdots,N},\]
where
\[\vt_n=\left[\cos\frac{2\pi n}{N},\sin\frac{2\pi n}{N}\right].\]
Then, for any searching point $\ms$, the traditional indicator function $\mathfrak{F}_{\mathrm{DSM}}(\ms)$ introduced in \cite{IJZ1,IJZ2,LZ} is given by
\[\mathfrak{F}_{\mathrm{DSM}}(\ms):=\frac{|\langle\psi_{\infty}(\md,\vt_n),\exp(-ik\vt_n\cdot\ms)\rangle_{L^2(\mathbb{S}^1)}|}{||\psi_{\infty}(\md,\vt_n)||_{L^2(\mathbb{S}^1)}||\exp(-ik\vt_n\cdot\ms)||_{L^2(\mathbb{S}^1)}},\]
where
\[\langle f(\vt_n),g(\vt_n)\rangle_{L^2(\mathbb{S}^1)}:=\sum_{n=1}^{N}f(\vt_n)\overline{g(\vt_n)}\quad\mbox{and}\quad||f(\vt_n)||_{L^2(\mathbb{S}^1)}:=\langle f(\vt_n),f(\vt_n)\rangle_{L^2(\mathbb{S}^1)}.\]
The locations of the inhomogeneities can be identified via the above indicator function when their permittivities differ from the background. In the permeability-contrast case, two peaks of large magnitude will appear in the neighborhood of each inhomogeneity instead of at their true locations. This is explained theoretically below.

\begin{thm}\label{TheoremDSM}For sufficiently large $N$, $\mathfrak{F}_{\mathrm{DSM}}(\ms)$ can be represented as follows:
\begin{equation}\label{Structure}
\mathfrak{F}_{\mathrm{DSM}}(\ms)=\frac{|\Psi(\ms)|}{\displaystyle\max_{\ms\in\mathbb{R}^2}|\Psi(\ms)|},
\end{equation}
where
\begin{equation}\label{StructurePsi}
\Psi(\ms):=\sum_{m=1}^{M}(r_m)^2\left(\frac{\mu_0}{\mu_m+\mu_0}\right)\left(\frac{\mx_m-\ms}{|\mx_m-\ms|}\cdot\md\right)J_1(k|\mx_m-\ms|).
\end{equation}
Here, $J_1$ denotes the Bessel function of order one of the first kind.
\end{thm}
\begin{proof}
Since $N$ is sufficiently large, applying asymptotic expansion formula \eqref{asymptotic} yields 
\begin{align*}
\Psi(\ms)&=\langle\psi_{\infty}(\md,\vt_n),\exp(-ik\vt_n\cdot\ms)\rangle_{L^2(\mathbb{S}^1)}=\sum_{n=1}^{N}\psi_{\infty}(\md,\vt_n)\exp(ik\vt_n\cdot\ms)\\
&=-\frac{k^2(1+i)}{4\sqrt{k\pi}}\sum_{n=1}^{N}\left(\sum_{m=1}^{M}(r_m)^2|\mB_m|(\md\cdot\mathbb{M}(\mx_m)\cdot\vt_n)\exp(ik\md\cdot\mx_m)\exp(-ik\vt_n\cdot\mx_m)\exp(ik\vt_n\cdot\ms)\right)\\
&=-\frac{k^2(1+i)}{4\sqrt{k\pi}}\sum_{m=1}^{M}(r_m)^2|\mB_m|\left(\frac{\mu_0}{\mu_m+\mu_0}\right)\exp(ik\md\cdot\mx_m)\sum_{n=1}^{N}(\md\cdot\vt_n)\exp(-ik\vt_n\cdot(\ms-\mx_m)).
\end{align*}
Notice that for $\mx\in\mathbb{R}^2$ and $\vx,\vt_n,\vt\in\mathbb{S}^1$, the following relation holds (see \cite{P-SUB3})
\[\sum_{n=1}^{N}(\vx\cdot\vt_n)\exp(-ik\vt_n\cdot\mx)=2\pi\int_{\mathbb{S}^1}(\vx\cdot\vt)\exp(-ik\vt\cdot\mx)d\vt=2\pi i\left(\vx\cdot\frac{\mx}{|\mx|}\right)J_1(k|\mx|),\]
we can immediately obtain
\[\Psi(\ms)=\frac{k^2\sqrt{\pi}(1-i)}{2}\sum_{m=1}^{M}(r_m)^2|\mB_m|\left(\frac{\mu_0}{\mu_m+\mu_0}\right)\exp(ik\md\cdot\mx_m)\left(\frac{\mx_m-\ms}{|\mx_m-\ms|}\cdot\md\right)J_1(k|\mx_m-\ms|).\]
Finally, since $|\mB_m|\equiv\pi$ and $|\exp(ik\md\cdot\mx_m)|\equiv1$, applying H{\"o}lder's inequality
\[|\langle\psi_{\infty}(\md,\vt_n),\exp(-ik\vt_n\cdot\ms)\rangle_{L^2(\mathbb{S}^1)}|\leq||\psi_{\infty}(\md,\vt_n)||_{L^2(\mathbb{S}^1)}||\exp(-ik\vt_n\cdot\ms)||_{L^2(\mathbb{S}^1)},\]
we can obtain \eqref{StructurePsi}. This completes the proof.
\end{proof}

\begin{rem}\label{Remark}
Based on the result in Theorem \ref{TheoremDSM}, we can observe that unlike the traditional result, $\mathfrak{F}_{\mathrm{DSM}}(\ms)$ should equal to $0$ (or small values) at the $\ms=\mx_m\in\Sigma_m$ because $J_1(0)=0$. Since $J_1(|x|)$ has its maximum value at $x=\pm1.8412$, the map of $\mathfrak{F}_{\mathrm{DSM}}(\ms)$ has two peaks of magnitude $1$ at the locations $\ms$ satisfy
\begin{equation}\label{Locations}
k|\mx_m-\ms|=1.8412\quad\mbox{and}\quad\frac{\mx_m-\ms}{|\mx_m-\ms|}=\pm\md.
\end{equation}
This means that location identification is highly depending on the value of $k$ and the direction of propagation $\md$, refer to Example \ref{Ex1}. Notice that, if the value of permeability $\mu_m$ is significantly larger than the others, the location of $\mx_m$ cannot be detected through the direct sampling method. This is the theoretical reason why the true locations of $\mx_m$ cannot be detected through the traditional direct sampling method.
\end{rem}

\section{Simulation results}\label{sec:4}
In this section, we present the results of some numerical simulations to support Theorem \ref{TheoremDSM}. For this purpose, we choose a set of three different small inhomogeneities $\Sigma_m$, $m=1,2,3$, with radii $r_m\equiv0.1$. Locations $\mx_m$ of $\Sigma_m$ are selected as $\mx_1=[0.7,0.5]$, $\mx_2=[-0.7,0.0]$, $\mx_3=[0.2,-0.5]$, the wavelength is chosen as $\lambda=0.4$, and the  direction of propagation is selected as $\md=[\cos(\pi/4),\sin(\pi/4)]$. It is worth mentioning that all far-field pattern data $\psi_{\infty}(\md,\vt_n)$ are generated by solving Foldy-Lax framework as presented in \cite{HSZ4}.

\begin{figure}[h]
\begin{center}
\includegraphics[width=0.495\textwidth]{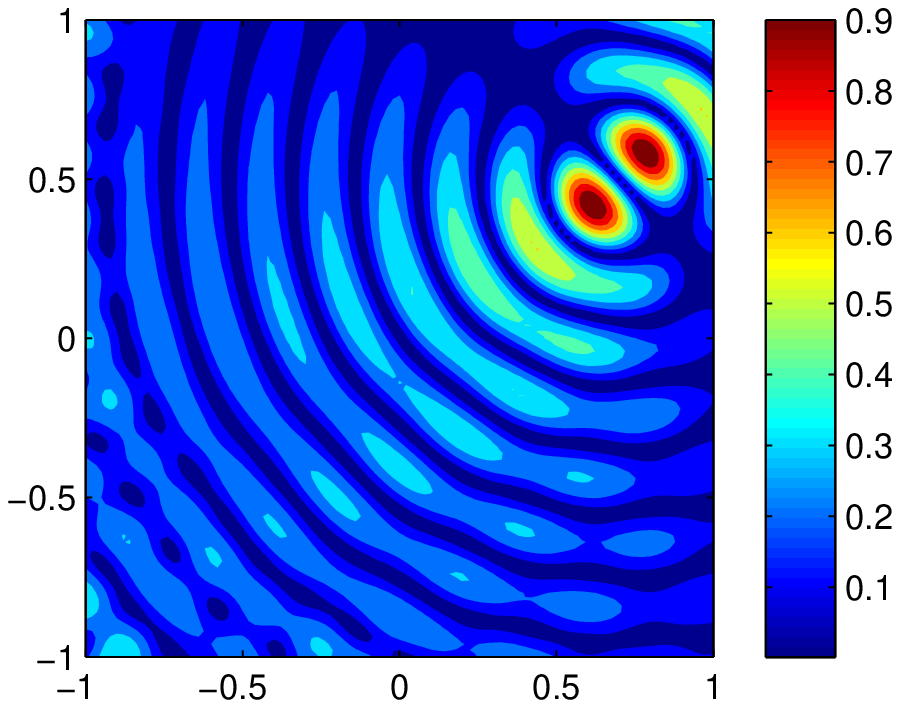}
\includegraphics[width=0.495\textwidth]{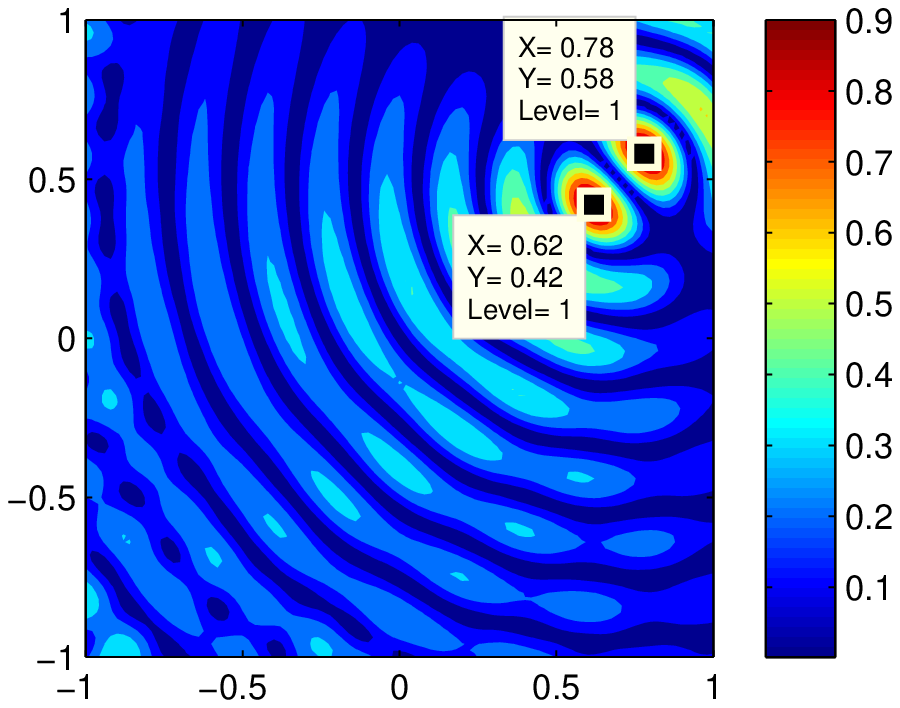}
\caption{\label{Result1}(Example \ref{Ex1}) Map of $\mathfrak{F}_{\mathrm{DSM}}(\mx)$ when $\lambda=0.4$.}
\end{center}
\end{figure}

\begin{ex}[Imaging of single inhomogeneity]\label{Ex1}
Figure \ref{Result1} shows map of $\mathfrak{F}_{\mathrm{DSM}}(\mx)$ when there is only one inhomogeneity $\Sigma_1$ whose permeability is $\mu_1=5$. As we discussed in Remark \ref{Remark}, instead of the true location $\mx_1$, the two peaks of largest magnitude appear at locations satisfying \eqref{Locations}. For example, the horizontal and vertical positions of these peaks are
\[0.7\pm\frac{1.8412}{k}\cos\frac{\pi}{4}=\set{0.6171,0.7829}\quad\mbox{and}\quad0.5\pm\frac{1.8412}{k}\sin\frac{\pi}{4}=\set{0.4171,0.5829},\]
respectively.
\end{ex}

\begin{figure}[h]
\begin{center}
\includegraphics[width=0.495\textwidth]{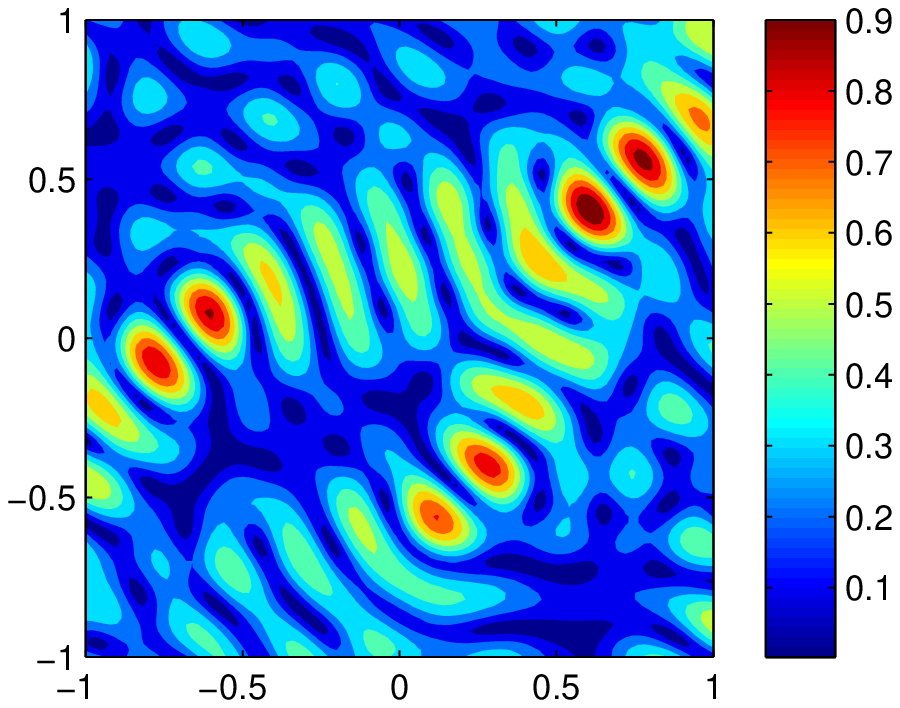}
\includegraphics[width=0.495\textwidth]{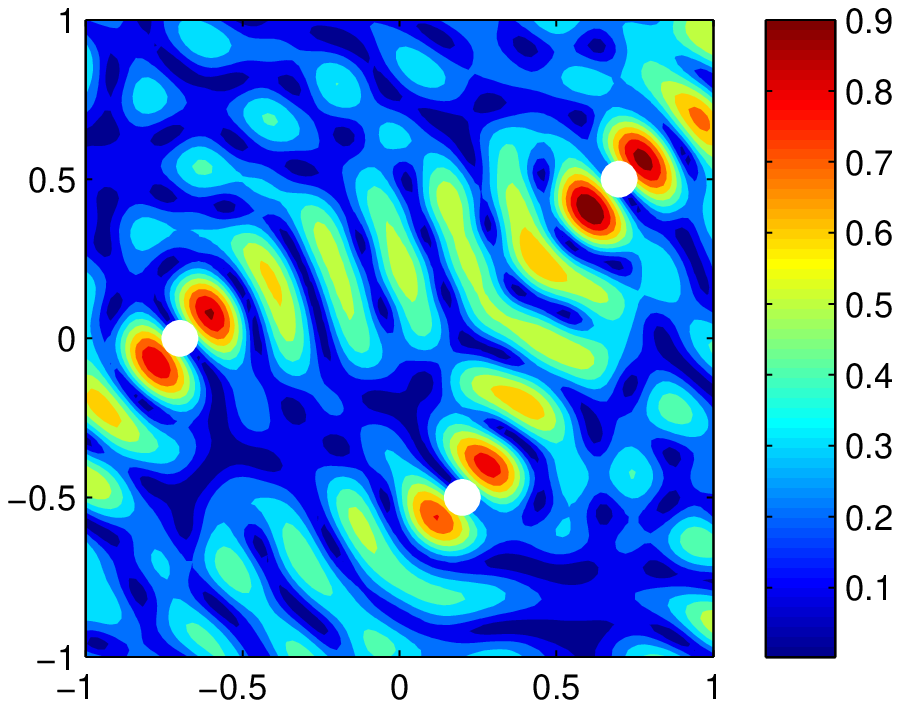}
\caption{\label{Result2}(Example \ref{Ex2}) Map of $\mathfrak{F}_{\mathrm{DSM}}(\mx)$ when $\lambda=0.40$. White circles denote the locations of $\mx_m$, $m=1,2,3$.}
\end{center}
\end{figure}

\begin{ex}[Imaging of multiple inhomogeneities with same value of permeability]\label{Ex2}
Figure \ref{Result2} shows map of $\mathfrak{F}_{\mathrm{DSM}}(\mx)$ with three inhomogeneities $\Sigma_m$ whose permeabilities are the same as $\mu_m\equiv5$, $m=1,2,3$. Similar to the result in Figure \ref{Result1}, the two peaks of largest magnitude appear away from the true locations $\mx_m$. Unfortunately, the appearance of large numbers amounts of artifacts disturbs the identification of the locations of $\Sigma_m$.
\end{ex}

\begin{figure}[h]
\begin{center}
\includegraphics[width=0.495\textwidth]{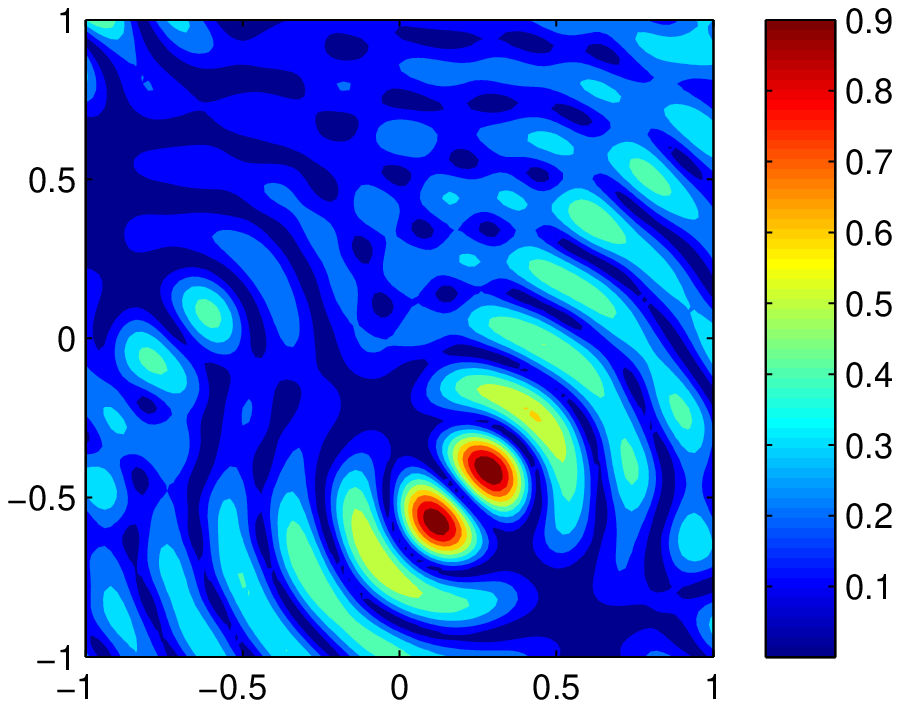}
\includegraphics[width=0.495\textwidth]{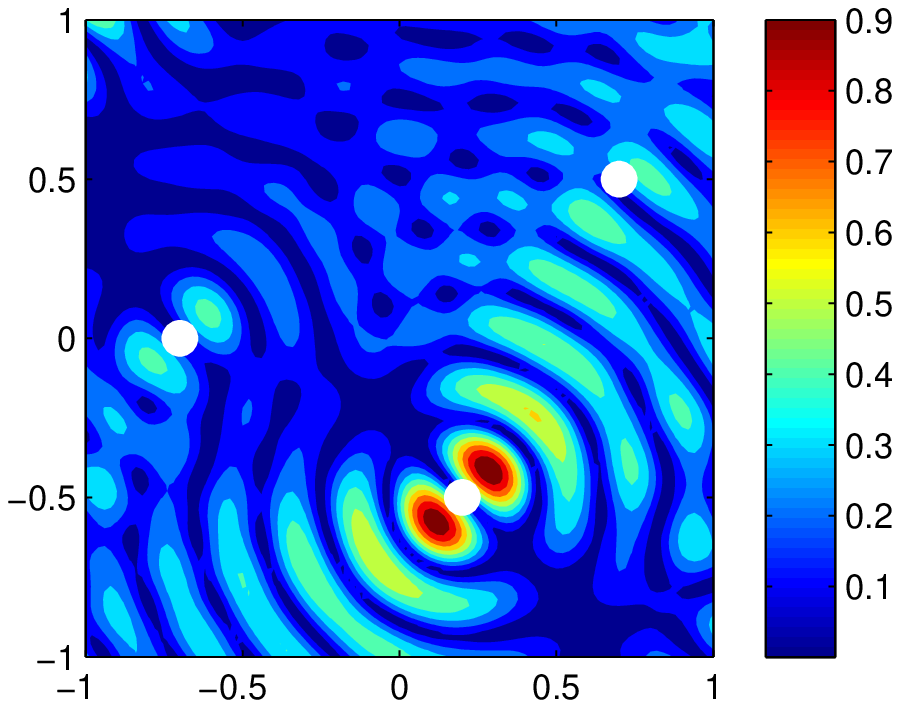}
\caption{\label{Result3}(Example \ref{Ex3}) Map of $\mathfrak{F}_{\mathrm{DSM}}(\mx)$ when $\lambda=0.40$. White circles denote the locations of $\mx_m$, $m=1,2,3$.}
\end{center}
\end{figure}

\begin{ex}[Imaging of multiple inhomogeneities with different value of permeabilities]\label{Ex3}
Figure \ref{Result3} shows map of $\mathfrak{F}_{\mathrm{DSM}}(\mx)$ for a case with three inhomogeneities $\Sigma_m$ whose permeabilities are $\mu_1=10$, $\mu_2=6$, and $\mu_3=2$. As we discussed in Remark \ref{Remark}, the existence of $\mx_3$ can be successfully identified because the value of permeability of $\Sigma_3$ is smaller than those of $\Sigma_1$ and $\Sigma_2$. However, owing to the existence of a large number of artifacts, it is very difficult to identify the existence of $\mx_2$. Furthermore, because the permeability of $\Sigma_1$ is very large, it is impossible to recognize the existence of $\Sigma_1$ through the direct sampling method.  
\end{ex}

\section{Conclusion}\label{sec:5}
We considered the direct sampling method for imaging small inhomogeneities when their permeabilities differ from the background. Throughout careful analysis, we established mathematical structure of indicator function by the Bessel function or order one, direction of propagation, and characteristic of inhomogeneity (size and permeability). Based on this, we observed that although the exact locations of such inhomogeneities cannot be detected via a traditional direct sampling method, they are located between two peaks of large magnitude. Thus, improving the direct sampling method to retrieve the exact locations of inhomogeneities in TE polarization will be an interesting research subject. Finally, extending current research to the real-world microwave imaging problem \cite{PKLS} will be the forthcoming work.

\section*{Acknowledgement}
The author would like to acknowledge two anonymous referees for their precious comments. This research was supported by the Basic Science Research Program of the National Research Foundation of Korea (NRF) funded by the Ministry of Education (No. NRF-2017R1D1A1A09000547) and the research program of Kookmin University in Korea.

\bibliographystyle{elsarticle-num-names}
\bibliography{References}
\end{document}